\documentclass[twoside,10pt]{article}

\usepackage{psfrag}
\usepackage{graphicx}
\usepackage{amssymb}
\usepackage{latexsym}
\usepackage{amsmath}

\newtheorem{theorem}{Theorem}
\newtheorem{proposition}{Proposition}
\newtheorem{definition}{Definition}
\newtheorem{lemma}{Lemma}

\newtheorem{remark}{Remark}
\newtheorem{example}{Example}

\date{}

\begin{document}

\title{A product construction for hyperbolic metric spaces}
\maketitle

\begin{center}
{\large Thomas Foertsch *\footnote{* Supported by the Deutsche Forschungsgemeinschaft (FO 353/1-1)} \hspace{1cm} Viktor Schroeder$^{\sharp}$\footnote{$\sharp$ 
Partially supported by the Suisse National Science Foundation}}
\footnote{2000 Mathematics Subject Classification. Primary 53C21} \\
\end{center}

\begin{abstract}
Given two pointed Gromov hyperbolic metric spaces $(X_i,d_i,z_i)$, $i=1,2$, and $\Delta \in \mathbb{R}^+_0$, we present a 
construction method, which yields another Gromov hyperbolic metric space $Y_{\Delta}=Y_{\Delta}((X_1,d_1,z_1),(X_2,d_2,z_2))$.
Moreover, it is shown that once $(X_i,d_i)$ is roughly geodesic, $i=1,2$, then there exists a ${\Delta}'\ge 0$ such that 
$Y_{\Delta}$ also is roughly geodesic for all $\Delta \ge {\Delta}'$.
\end{abstract}

\vspace{0.5cm}

%%%%%%%%%%%%%%%%%%%%%%%%%%%%%%%%%%%%%%%%%%%%%%%%%%%%%%%%%%%%%%%%%%%%%%%%%%%%%%%%%%%%%%%%%%%%%
%%%%%%%%%%%%%%%%%%%%%%%%%%%%%%%%%%%%%%%%%%%%%%%%%%%%%%%%%%%%%%%%%%%%%%%%%%%%%%%%%%%%%%%%%%%%%

\section{Introduction}

\label{sec-intro}

A metric space $(X,d)$ is called  $\delta$-hyperbolic, $\delta \ge 0$, if for all $x,y,z,w\in X$ it holds
\begin{equation} \label{eqn-def-hyperbolicity}
d(x,y)+d(z,w) \; \le \; \max \{ d(x,z)+d(y,w),d(x,w)+d(y,z)\} \; + \; 2\delta
\end{equation}
and said to be Gromov hyperbolic, if it is $\delta$-hyperbolic for some $\delta \ge 0$. \\

Let $(X_i,d_i)$ be Gromov hyperbolic metric spaces and fix $z_i\in X_i$, $i=1,2$. For $\Delta \ge 0$ consider the set
$Y_{\Delta}=Y_{\Delta}(X_1,d_1,z_1,X_2,d_2,z_2)$ defined via
\begin{displaymath}
Y_{\Delta} \; := \; \{ (x_1,x_2)\in X_1\times X_2 \; | \; |d_1(x_1,z_1)-d_2(x_2,z_2)|\le \Delta\} \subset X \; := \; X_1\times X_2.
\end{displaymath}
On $Y_{\Delta}$ we consider the metric $d_m|_{Y_{\Delta}\times Y_{\Delta}}$ which is the restriction of the $l_{\infty}$-product metric 
$d_m:X\times X\longrightarrow \mathbb{R}^+_0$
\begin{displaymath}
d_m\Big( (x_1,x_2),(x_1',x_2')\Big) \; := \; \max \{d_1(x_1,x_1'),d_2(x_2,x_2')\},
\end{displaymath}
for all $x_1,x_1'\in X_1, x_2,x_2'\in X_2$, to $Y_{\Delta}\times Y_{\Delta}\subset X\times X$. \\

Our paper is based on the following elementary observation which we refer to as
\begin{theorem} \label{theo-gen}
Let $(X_i,d_i)$ be Gromov hyperbolic metric spaces and $z_i\in X_i$, $i=1,2$. Then $(Y_{\Delta},d_m|_{Y_{\Delta}\times Y_{\Delta}})$,
as introduced above, is Gromov hyperbolic.
\end{theorem}

With Theorem \ref{theo-gen} at hand we further prove the 

\begin{theorem} \label{theo-r-geod}
Let $(X_i,d_i)$ be roughly geodesic, Gromov hyperbolic metric spaces, $z_i\in X_i$, $i=1,2$, then there exists $\tilde{\Delta}\ge 0$ such that
$(Y_{\Delta},d_m|_{Y_{\Delta}\times Y_{\Delta}})$ is roughly geodesic for all $\Delta \ge \tilde{\Delta}$ (and hyperbolic due
to Theorem \ref{theo-gen}). Moreover, its boundary  at infinity is naturally homeomorphic to 
${\partial}_{\infty}(X_1,d_1)\times {\partial}_{\infty}(X_2,d_2)$.
\end{theorem}

For precise definitions of the boundary at infinity and rough geodesics see Section \ref{sec-basics}.

\begin{remark}
\begin{description}
\item[(i)] Both theorems can be formulated for a finite number of factors.
\item[(ii)] For metric spaces $(X_i,d_i)$ with nonempty boundaries at infinity, the  Theorems 
\ref{theo-gen} and \ref{theo-r-geod} have analogues in the limit case, that the fixed points $z_i\in X_i$ converge at infinity. Those will
precisely be stated in Section \ref{sec-limit}.
\item[(iii)] An analogue of Theorem \ref{theo-r-geod} in the setting of geodesic metric spaces has been studied in \cite{fs2}. 
\end{description}
\end{remark}

{\bf Outline of the paper:} In Section \ref{sec-theo-gen} we prove Theorem \ref{theo-gen}. In Section \ref{sec-basics} we recall
some basic definitions and facts on Gromov hyperbolic metric spaces, which will be used in Section \ref{sec-theo-r-geod} when 
proving Theorem \ref{theo-r-geod} and Section \ref{sec-limit} when we state the above mentioned ``limit case analogues'' of the
Theorems \ref{theo-gen} and \ref{theo-r-geod}. \\

{\bf Acknowledgment:} It is a pleasure to thank Mario Bonk for useful discussions and valuable comments on an earlier version of this paper.

%%%%%%%%%%%%%%%%%%%%%%%%%%%%%%%%%%%%%%%%%%%%%%%%%%%%%%%%%%%%%%%%%%%%%%%%%%%%%%%%%%%%%%%%%%%%%%%%%%%%%%%%%%%%%%%%%%%%%%%%%%%%%%%%%%%%%%
%%%%%%%%%%%%%%%%%%%%%%%%%%%%%%%%%%%%%%%%%%%%%%%%%%%%%%%%%%%%%%%%%%%%%%%%%%%%%%%%%%%%%%%%%%%%%%%%%%%%%%%%%%%%%%%%%%%%%%%%%%%%%%%%%%%%%%

\section{The Proof of Theorem \ref{theo-gen}}

\label{sec-theo-gen}

{\bf Proof of Theorem \ref{theo-gen}:} First of all note that if for a metric space $(X,d)$ there exist $z\in X$ and $\tilde{\delta} \ge 0$
such that for all $x,y,w\in X$ it holds
\begin{equation} \label{eqn-4point->3point}
d(x,y) \; + \; d(w,z) \; \le \; 
\max \{ d(x,w)+d(y,z),d(x,z)+d(y,z)\} \; + \; 2\tilde{\delta}
\end{equation} 
then $(X,d)$ is $2\tilde{\delta}$-hyperbolic (see e.g. \cite{g}). \\

Let now $(X_i,d_i)$ be ${\delta}_i$-hyperbolic metric spaces and set $\delta := \max \{{\delta}_1,{\delta}_2\}$. In order to 
show that $(Y_{\Delta},d_m)$ is hyperbolic, we show that inequality (\ref{eqn-4point->3point}) holds for
$d=d_m$, $z=(z_1,z_2)\in Y_{\Delta}$, $\tilde{\delta}:=\delta +\frac{\Delta}{2}$ and $x,y,z\in Y_{\Delta}$ arbitrary: \\

Without loss of generality we assume $d_1(x_1,y_1)\ge d_2(x_2,y_2)$. Now, due to the definition of $Y_{\Delta}$ we find
\begin{displaymath}
d_m(w,z) \; = \; \max \{ d_1(w_1,z_1),d_2(w_2,z_2)\} \; \le \; d_1(w_1,z_1) \; + \; \Delta .
\end{displaymath}
Thus with the $\delta$-hyperbolicity of the first factor we get
\begin{eqnarray*}
& & d_m(x,y)+d_m(w,z) \\
& \le & d_1(x_1,y_1) \; + \; d_1(w_1,z_1) \; + \; \Delta \\
& \le & \max \{ d_m(x,w)+d_m(y,z),d_m(x,z)+d_m(y,z)\} \; + \; 2\delta \; + \; \Delta .
\end{eqnarray*}
\hfill $\Box$ \\

%%%%%%%%%%%%%%%%%%%%%%%%%%%%%%%%%%%%%%%%%%%%%%%%%%%%%%%%%%%%%%%%%%%%%%%%%%%%%%%%%%%%%%%%%%%%%%%%%%%%%%%%%%%%%%%%%%%%%%%%%%%%%%%%%%%%%%
%%%%%%%%%%%%%%%%%%%%%%%%%%%%%%%%%%%%%%%%%%%%%%%%%%%%%%%%%%%%%%%%%%%%%%%%%%%%%%%%%%%%%%%%%%%%%%%%%%%%%%%%%%%%%%%%%%%%%%%%%%%%%%%%%%%%%%

\section{Some basic definitions}

\label{sec-basics}

\subsection{Hyperbolicity and rough geodesics}

Let $(X_1,d_1),(X_2,d_2)$ be metric spaces. A map $f:X_1\longrightarrow X_2$ is called a quasi-isometric embedding, if 
there exist $\lambda \ge 1$, $k\ge 0$ such that for all $x,x'\in X_1$ it holds
\begin{displaymath}
\frac{1}{\lambda}d_1(x,x') \; - \; k \; \le \; d_2\Big( f(x),f(x')\Big) \; \le \; \lambda d_1(x,x') \; + \; k.
\end{displaymath}
If $k=0$, $f$ is called bilipschitz, while, for $\lambda =1$, $f$ is said to be a $k$-rough-isometric embedding. For $k=0$ and 
$\lambda =1$ the embedding is called isometric. \\

Let $f:X_1\longrightarrow X_2$ be a rough isometric embedding. Then, if $(X_2,d_2)$ is hyperbolic, so is $(X_1,d_1)$. In case
$d_1$ and $d_2$ are length metrics, then the same holds when replacing the rough-isometric embedding through a quasi-isometric
embedding (This non-trivial but by now standard result may be found in, for instance, \cite{brih} or \cite{bubui}). \\

A ($k$-rough) geodesic $\gamma$ in $(X,d)$ connecting $x\in X$ to $x'\in X$ is a ($k$-rough) isometric embedding
$\gamma :[\alpha ,\omega ]\longrightarrow X$ such that $\gamma(\alpha )=x$ and $\gamma (\omega )=x'$. The metric space
$(X,d)$ is called ($k$-rough) geodesic if for all $x,x'\in X$ there exists a ($k$-rough) geodesic connecting $x$ to $x'$.
If there exists a $k\ge 0$ such that $(X,d)$ is $k$-rough geodesic, then $(X,d)$ is said to be rough geodesic. \\

According to \cite{bos} a metric space $(X,d)$ is called $k$-almost geodesic, if for every $x,y\in X$ and every $t\in [0,d(x,y)]$,
there exists $w\in X$ such that $|d(x,w)-t|\le k$ and $|d(y,w)-(d(x,y)-t)|\le k$. $(X,d)$ is said to be almost geodesic, if 
there exists $k\ge 0$ such that it is $k$-almost geodesic. Note that a hyperbolic metric space $(X,d)$ is alomst geodesic if and
only if it is rough geodesic (compare Proposition 5.2 in \cite{bos}). \\

For geodesic metric spaces there are a number of equivalent characterizations of hyperbolicity using the 
geometry of geodesic triangles (see e.g. \cite{brih} and \cite{fs2}). All of those characterizations have analogues in the 
rough geodesic setting. Here we state the corresponding results, we are going to make use of in the following:

\begin{definition}
A $k$-roughly geodesic metric
space $(X,d)$ is said to be $(\delta ,k)$-hyperbolic if each side of any $k$-roughly geodesic triangle in $(X,d)$ is contained in 
the $\delta$-neighborhood of the union of the other two sides. 
\end{definition}
Just along the lines of the proof of the corresponding statement for geodesic spaces (see e.g. \cite{brih}) one proves the
\begin{proposition} \label{prop-delta-k-hyp}
Let $(X,d)$ be a $k$-roughly geodesic space. Then the following are equivalent:
\begin{description}
\item[(1)] $(X,d)$ is Gromov-hyperbolic.
\item[(2)] there exists a $\delta \in \mathbb{R}_0^+$ such that $(X,d)$ is $(\delta ,k)$-hyperbolic.
\end{description}
\end{proposition}

Let $X$ be a metric space and $x,y,z\in X$. Then there exist unique $a,b,c\in \mathbb{R}^+_0$ such that 
\begin{displaymath}
d(x,y)=a+b, \hspace{0.5cm} d(x,z)=a+c \hspace{0.5cm} \mbox{and} \hspace{0.5cm} d(y,z)=b+c . 
\end{displaymath}
In fact those numbers are given through
\begin{displaymath}
a \; = \; (y\cdot z)_x \; , \hspace{0.5cm} b \; = \; (x\cdot z)_y   \hspace{0.5cm}
\mbox{and} \hspace{0.5cm} c \; = \; (x\cdot y)_z \; , 
\end{displaymath}
where, for instance, 
\begin{displaymath}
(y\cdot z)_x \; = \; \frac{1}{2} \Big[ d(y,x) \; + \; d(z,x) \; - \; d(y,z) \Big] .
\end{displaymath}
In the case that $X$ is $k$-roughly geodesic we may consider a $k$-roughly geodesic triangle 
$\overline{xy} \cup \overline{xz} \cup \overline{yz} \subset X$,
where for example $\overline{xy}$ denotes a $k$-roughly geodesic segment connecting $x$ to $y$. Given such a triangle we 
write $\tilde{x}:={\gamma}_{yz}(b)$, $\tilde{y}:={\gamma}_{xz}(a)$ and $\tilde{z}:={\gamma}_{xy}(a)$. Note that
for geodesic triangles it holds ${\gamma}_{xz}(a)={\gamma}_{xz}^{-1}(c)$. In the case that $(X,d)$ is only 
$k$-roughly geodesic we still have e.g. $d({\gamma}_{xz}(a),{\gamma}_{xz}^{-1}(c))\le 2k$. \\

Similar to Lemma 1 i) in \cite{fs2} one proves the
\begin{lemma} \label{lemma-delta-k-hyperbolic}
If $(X,d)$ is $(\delta ,k)$-hyperbolic, then
\begin{displaymath}
d(z,\tilde{z}) \; \le c \; + \; 2\delta \; + \; 4k, \hspace{1cm} 
d\Big( {\gamma}_{xy}(t),{\gamma}_{xz}(t)\Big) \; \le \; 4\delta \; + \; 15k \;\;\; \forall \; t\in [0,a]
\end{displaymath}
and the points $\tilde{x}$, $\tilde{y}$ and $\tilde{z}$ have pairwise distance $\le 4\delta +15k$.
\end{lemma}

%%%%%%%%%%%%%%%%%%%%%%%%%%%%%%%%%%%%%%%%%%%%%%%%%%%%%%%%%%%%%%%%%%%%%%%%%%%%%%%%%%%%%%%%%%%%%

\subsection{The boundary at infinity and Busemann functions}

Given a hyperbolic space $(X,d)$ there are various ways to attach a boundary at infinity ${\partial}_{\infty}X$ to $X$.
In this paper we define ${\partial}_{\infty}X$ in the following way: \\
We choose a basepoint $z\in X$ and say that a sequence $\{ x^i{\}}_{i\in \mathbb{N}}$ of points in $X$ 
{\it converges to infinity}, if 
\begin{displaymath}
\liminf\limits_{i,j\longrightarrow \infty} \; (x^i\cdot x^j)_z \; = \; \infty .
\end{displaymath}
Two sequences  $\{ x^i{\}}_{i\in \mathbb{N}}$ and  $\{ y^i{\}}_{i\in \mathbb{N}}$ converging to infinity
are equivalent, \linebreak  $\{ x^k{\}}_{k\in \mathbb{N}} \sim \{ y^k{\}}_{k\in \mathbb{N}}$ if
\begin{displaymath}
\liminf\limits_{i,j\longrightarrow \infty} \; (x^i\cdot y^j)_z \; = \; \infty .
\end{displaymath}
One shows that $\sim$ is an equivalence relation and defines $\partial X$ as the set of equivalence classes.
We write $[\{ x^i\} ]\in \partial X$ for the corresponding class. \\ 

For $v\in \partial X$ and $r>0$ one defines
\begin{displaymath}
U(v,r) \; := \; 
\Big\{ {\textstyle
w\in \partial X \; \Big| \; \exists \{ x^k\} , \{ y^k\} \; \mbox{s.t.} \;
[\{ x^k\} ]=v, \;  [\{ y^k\} ]=w, \; \liminf\limits_{k,l\longrightarrow \infty} \; (x^k\cdot y^l)_z \; > \; r }
\Big\} . 
\end{displaymath} 
On $\partial X$ we consider the topology generated by $U(v,r)$, $v\in \partial X$, $r>0$. \\

Let now $(X,d)$ be $k$-roughly geodesic, then there exists a $k'=k'(k,\delta )$ with $\delta$ as in 
equation (\ref{eqn-def-hyperbolicity}) such that for every $x\in X$ there exists a $k'$-rough geodesic 
${\gamma}_{xu}: [0,\infty )\longrightarrow X$ with ${\gamma}_{xu}(0)=x$ and $[\{ {\gamma}_{xy}(i)\} ]=u$ 
(see \cite{bos}). Such rays are said to connect $x$ to $u$. \\
We now fix such a $k'$-roughly geodesic ray ${\gamma}_{zu}$ connecting $z\in X$ to $u\in {\partial}_{\infty}X$ and define
the Busemann function $B_{{\gamma}_{zu}} : X\longrightarrow \mathbb{R}$ associated to the ray ${\gamma}_{zu}$ via
\begin{displaymath}
B_{{\gamma}_{zu}} (x) \; := \; \liminf\limits_{t\longrightarrow \infty} 
\Big[ d\Big( x,{\gamma}_{zu}(t) \Big) \; - \; t \Big] .
\end{displaymath}
Note that the limit inferior always exists, while the limit itself necessarily only exists once ${\gamma}_{zu}$
is a geodesic.

%%%%%%%%%%%%%%%%%%%%%%%%%%%%%%%%%%%%%%%%%%%%%%%%%%%%%%%%%%%%%%%%%%%%%%%%%%%%%%%%%%%%%%%%%%%%%%%%%%%%%%%%%%%%%%%%%%%%%%%%%%%%%%%%%%%%%%
%%%%%%%%%%%%%%%%%%%%%%%%%%%%%%%%%%%%%%%%%%%%%%%%%%%%%%%%%%%%%%%%%%%%%%%%%%%%%%%%%%%%%%%%%%%%%%%%%%%%%%%%%%%%%%%%%%%%%%%%%%%%%%%%%%%%%%

\section{The proof of Theorem \ref{theo-r-geod}}

\label{sec-theo-r-geod}

In this section we provide the 

{\bf Proof of Theorem \ref{theo-r-geod}:}
Let $(X_i,d_i)$ be $k_i$-roughly geodesic, ${\delta}_i$-hyperbolic, $i=1,2$, and set $k:=\max \{ k_1,k_2\}$ as well as 
$\delta := \max \{ {\delta}_1,{\delta}_2\}$. We show that for 
$\Delta \ge 4k$ the space $(Y_{\Delta},d_m)$ is $K(\Delta ,k,\delta )$-almost geodesic and therefore roughly geodesic due to 
Theorem \ref{theo-gen} and Proposition 5.2 in \cite{bos}. \\
Thus we have to show that for $\Delta \ge 4k$ there exists $K\ge 0$ such that for all $x,y\in Y_{\Delta}$, $t\in [0,d_m(x,y)]$
there exists $w\in Y_{\Delta}$ such that
\begin{displaymath}
\begin{array}{rcccl}
t \; - \; K & \le & d_m(x,w) & \le & t \; + \; K \hspace{0.5cm} \mbox{and} \hspace{1cm} (*)\\
& & & & \\
d_m(x,y)-t-K & \le & d_m(y,w) & \le & d_m(x,y)-t-K . 
\end{array}
\end{displaymath}
(1) W.l.o.g. we assume $d_m(x,y)=d_1(x_1,y_1)\ge d_2(x_2,y_2)$. \\
(2) W.l.o.g. we assume $\Delta < t < d_m(x,y)-\Delta$. This can be done, since for $\Delta <t$ we may set $w:=x$ while for
$t> d_m(x,y) - \Delta $ we may set $w:=y$ and the inequalities above trivally hold once $K\ge 2\Delta$. \\
(3) Now set
\begin{eqnarray*}
a_i & := & \frac{1}{2} \Big( d_i(x_i,y_i) \; + \; d_i(x_i,z_i) \; - \; d_i(y_i,z_i)\Big) , \\
b_i & := & \frac{1}{2} \Big( d_i(y_i,x_i) \; + \; d_i(y_i,z_i) \; - \; d_i(x_i,z_i)\Big) , \\
c_i & := & \frac{1}{2} \Big( d_i(z_i,x_i) \; + \; d_i(z_i,y_i) \; - \; d_i(x_i,y_i)\Big) , \\
\end{eqnarray*}
$i=1,2$. W.l.o.g. we assume $\Delta < t \le a_1$,
let ${\gamma}_{z_ix_i}:[0,d_i(z_i,x_i)]\longrightarrow X_i$ be $2k$-rough geodesics connecting $x_i$ to $z_i$, $i=1,2$ and set
\begin{displaymath}
w \; := \; \Big( {\gamma}_{z_1x_1}(d_1(x_1,z_1)-t), {\gamma}_{z_2x_2}(d_1(x_1,z_1)-t)\Big) .
\end{displaymath}
Note that since $\Delta \ge 4k$ it follows that $w\in Y_{\Delta}$. \\
Moreover, from the definition of $w$ it is clear that we have
\begin{equation} \label{eqn-w2-w1}
d_2(w_2,x_2) \; \le \; d_1(w_1,x_1) \; + \; \Delta \; + \; 4k.
\end{equation}
$(X_i,d_i)$ is hyperbolic, $i=1,2$. Thus, due to Proposition \ref{prop-delta-k-hyp}, there exists $\tilde{\delta}_i$ such that
$(X_i,d_i)$ is $(\tilde{\delta}_i,2k)$-hyperbolic. Setting $\tilde{\delta} := \max \{ \tilde{\delta}_1,\tilde{\delta}_2\}$ and
${\delta}':=4\tilde{\delta}+30k$ yields, due to Lemma \ref{lemma-delta-k-hyperbolic}, 
\begin{displaymath}
\begin{array}{rcccl}
t-2k-{\delta}' & \le & d_1(w_1,x_1) & \le & t+2k+{\delta}' \hspace{0.5cm} \mbox{and} \hspace{1cm} (**)\\
d_1(x_1,y_1)-t-2k-{\delta}' & \le & d_1(w_1,y_1) & \le & d_1(x_1,y_1)-t+2k-{\delta}' . 
\end{array}
\end{displaymath}
Now we consider the following two cases: \\
(i) $d_2(x_2,z_2)-[d_1(x_1,z_1)-t]\le a_2$: In this case Lemma \ref{lemma-delta-k-hyperbolic} yields
\begin{displaymath}
\begin{array}{rcccl}
t-2k-{\delta}' & \le & d_2(w_2,x_2) & \le & t+2k+{\delta}' \hspace{0.5cm} \mbox{and} \\
d_2(x_2,y_2)-t-2k-{\delta}' & \le & d_2(w_2,y_2) & \le & d_2(x_2,y_2)-t+2k-{\delta}' . 
\end{array}
\end{displaymath}
(ii) $d_2(x_2,z_2)-[d_1(x_1,z_1)-t]> a_2$: Of course we have $d_1(x_1,z_1)-t\ge d_1(x_1,z_1)-a_1=c_1$, hence
\begin{eqnarray*}
d_2(y_2,z_2) \; - \; [d_1(x_1,z_1)-t] & \le & d_2(y_2,z_2) -c_1 \\
& \le & d_1(y_1,z_1)+\Delta -c_1 \; = \; b_1+\Delta .
\end{eqnarray*}
Thus, due to  $d_2(x_2,z_2)-[d_1(x_1,z_1)-t]> a_2$ and Lemma \ref{lemma-delta-k-hyperbolic} we conclude
\begin{displaymath}
d_2(w_2,y_2) \; \le \; b_1 \; + \Delta \; + \; {\delta}' .
\end{displaymath}
From (i) and (ii) and the inequalities $(**)$ as well as (\ref{eqn-w2-w1}) it follows that the inequalities $(*)$ hold for
a $K\ge 0$ sufficiently large. \\

For the part of the proof concerning the boundary at infinity, we refer the reader to \cite{fs2}, where the geodesic case is treated.
\hfill $\Box$

%%%%%%%%%%%%%%%%%%%%%%%%%%%%%%%%%%%%%%%%%%%%%%%%%%%%%%%%%%%%%%%%%%%%%%%%%%%%%%%%%%%%%%%%%%%%%%%%%%%%%%%%%%%%%%%%%%%%%%%%%%%%%%%%%%%%%%
%%%%%%%%%%%%%%%%%%%%%%%%%%%%%%%%%%%%%%%%%%%%%%%%%%%%%%%%%%%%%%%%%%%%%%%%%%%%%%%%%%%%%%%%%%%%%%%%%%%%%%%%%%%%%%%%%%%%%%%%%%%%%%%%%%%%%%

\section{The limit case and final remarks}

\label{sec-limit}

In this section we state the Theorems \ref{theo-gen-lim} and \ref{theo-r-geod-lim} corresponding to the Theorems
\ref{theo-gen} and \ref{theo-r-geod} when fixing points in the boundary at infinity of the factors rather than points in the 
interior. \\

Let therefore $(X_1,d_1)$ and $(X_2,d_2)$ be roughly geodesic hyperbolic metric spaces with non-empty boundaries at infinity.
Fix $u_i\in {\partial}_{\infty}X_i$ as well as Busemann functions $B_i:X_i\longrightarrow \mathbb{R}$, associated to 
roughly geodesic rays ${\gamma}_i$ converging to $u_i$, $i=1,2$. This time we consider the sets
\begin{displaymath}
Y_{\Delta} \; := \; \{ (x_1,x_2)\in X_1\times X_2 \; | \; |B_1(x_1)-B_2(x_2)|\le \Delta \}, \hspace{0.5cm} \Delta \ge 0.
\end{displaymath}

With this notation the following theorems hold:

\begin{theorem} \label{theo-gen-lim}
Let $(X_i,d_i)$ be Gromov hyperbolic metric spaces with nonempty boundaries at infinity such that two Busemann functions
$B_i:X_i\longrightarrow \mathbb{R}$ associated to rough geodesic rays are defined. Then, for all $\Delta \ge 0$, $(Y_{\Delta},d_m|_{Y_{\Delta}\times Y_{\Delta}})$
also is hyperbolic.
\end{theorem}

\begin{theorem} \label{theo-r-geod-lim}
Let $(X_i,d_i)$ be roughly geodesic, Gromov hyperbolic metric spaces and $B_i:X_i\longrightarrow \mathbb{R}$ Busemann functions on
$X_i$, $i=1,2$. Then there exists $\tilde{\Delta} \ge 0$ such that $(Y_{\Delta},d_m|_{Y_{\Delta}\times Y_{\Delta}})$ also
is roughly geodesic and hyperbolic for all $\Delta \ge \tilde{\Delta}$. Moreover, ${\partial}_{\infty}(Y_{\Delta},d_m|_{Y_{\Delta}\times Y_{\Delta}})$ is 
naturally homeomorphic to the smashed product ${\partial}_{\infty}(X_1,d_1)\wedge {\partial}_{\infty}(X_2,d_2)$.
\end{theorem}

\begin{remark}
The smashed product $\wedge$ is a standard construction 
for pointed topological spaces (see e.g. \cite{m}). Let $(U_1,u_1)$, $(U_2,u_2)$ be two pointed spaces then the smashed product 
$U_1\wedge U_2$ is defined as $U_1\times U_2/U_1\vee U_2$, where $U_1\times U_2$ 
is the usual product and 
\begin{displaymath}
U_1\vee U_2 \; = \; \Big( \{ u_1\} \times U_2\Big) \; \cup \Big( U_1 \times \{ u_2\} \Big) \; \subset \; U_1\times U_2
\end{displaymath}
is the wedge product canonically embedded in $U_1\times U_2$. Thus $U_1\wedge U_2$ is obtained from $U_1\times U_2$
by collapsing $U_1\vee U_2$ to a point. For example $S^m\wedge S^n=S^{m+n}$.  
\end{remark}

The proofs of these theorems go just along the lines of the proofs of the corresponding Theorems \ref{theo-gen} and \ref{theo-r-geod} 
when fixing points in the interior rather than the boundary. For the part of Theorem \ref{theo-r-geod} concerning the boundary at infinity
we refer the reader to \cite{fs2}, where the analogue in the geodesic setting is proved. \\

We finally point out that when starting off with two proper geodesic metric spaces one has to consider the length metric $d$ induced by $d_m$
on $Y_0$, in order to obtain a proper geodesic space again. In this case, we might as well endow $Y_0$ with the length metric induced by the Euclidean 
product metric $d_e$ instead of the maximum metric $d_m$. Since both are length spaces which are bilipschitz related, one of them is 
Gromov hyperbolic if and only if the other one is. \\
In fact, when starting off with two Riemannian manifolds and fixing points at infinity, the construction using the Euclidean product metric
has the advantage that it once again yields a Riemannian manifold (compare e.g. \cite{fs1}). However, we emphazise that for neither of
the Theorems \ref{theo-gen}, \ref{theo-r-geod}, \ref{theo-gen-lim} and \ref{theo-r-geod-lim} we might replace the maximum metric through the 
Euclidean metric. This is, for instance, seen in the 

\begin{example} \label{example}
Consider two copies of the real hyperbolic space $\mathbb{H}^2$. Fix points $u_i\in {\partial}_{\infty}\mathbb{H}^2$, Busemann functions
$B_i$ associated to geodesic rays ${\gamma}_i$ converging to $u_i$, $i=1,2$, and 
consider sequences of points $\{x^n=(x^n_1,x^n_2)\}$, $\{y^n=(y^n_1,y^n_2)\}$, $\{z^n=(z^n_1,z^n_2)\}$ and $\{w^n=(w^n_1,w^n_2)\}$ 
such that $x^n_1=x^n_2$, $y^n_1=y^n_2$, $z^n_1=z^n_2$, $B_i(z^n_i)=B_i(y^n_i)$, $d_i(x^n_i,y^n_i)=d_i(x^n_i,z^n_i)=\frac{1}{2}d_i(y_i^n,z_i^n)$, $w^n_1=y^n_1$ and 
$w^n_2=z^n_2$ for all $n\in \mathbb{N}$, $i=1,2$ as well as $d_i(y^n_i,z^n_i) \stackrel{n\rightarrow \infty}{\longrightarrow} \infty$, $i=1,2$. \\
We claim that $(Y_0(\mathbb{H}^2,u_i,B_i),d_e)$ is not hyperbolic. Suppose the contrary, then there exists a $\delta \ge 0$ such that for all 
$n\in \mathbb{N}$
\begin{displaymath}
\begin{array}{crcl}
 & & &  d_e(y^n,z^n) \; + \; d_e(x^n,w^n) \\ 
& & & \\
& & \le & \max \{ d_e(x^n,y^n)+d(z^n,w^n),d_e(y^n,w^n)+d(x^n,z^n)\} \; + \; 2\delta \\
& & & \\
\Longleftrightarrow & d_e(y^n,z^n) & \le & \max \{ d_e(z^n,w^n),d(y^n,w^n)\} \; + \; 2\delta \\
& & & \\
\Longleftrightarrow & \sqrt{2} \, d_1(y^n_1,z^n_1) & \le & d_1(y^n_1,z^n_1) \; + \; 2\delta  ,
\end{array}
\end{displaymath}
which contradicts our choices of sequences. 
\end{example}

%%%%%%%%%%%%%%%%%%%%%%%%%%%%%%%%%%%%%%%%%%%%%%%%%%%%%%%%%%%%%%%%%%%%%%%%%%%%%%%%%%%%%%%%%%%%%%%%%%%%%%%%%%%%%%%%%%%%%%%%%%%%%%%%%%%%%%
%%%%%%%%%%%%%%%%%%%%%%%%%%%%%%%%%%%%%%%%%%%%%%%%%%%%%%%%%%%%%%%%%%%%%%%%%%%%%%%%%%%%%%%%%%%%%%%%%%%%%%%%%%%%%%%%%%%%%%%%%%%%%%%%%%%%%%


\begin{thebibliography}{cccccccc}
%\bibitem[BraF]{braf} N. Brady \& B. Farb, {\it Filling-Invariants at Infinity for Manifolds of Nonpositive Curvature},
%Trans. Amer. Math. Soc., Vol. 350, Num. 8, 3393-3405, 1998
\bibitem[BriH]{brih}  M. Bridson \& A. Haefliger, {\it Metric spaces of non-positive curvature}, Springer Verlag Berlin 1999
\bibitem[BoS]{bos} M. Bonk \& O. Schramm,  {\it Embeddings of Gromov hyperbolic spaces},
Geom. Funct. Anal. 10, No.2, 266-306, 2000
\bibitem[BuBuI]{bubui} D. Burago \& Y. Burago \& S. Ivanov, {\it A course in Metric Geometry}, Graduate Studies in Mathematics, 
Vol. 33, Amer. Math. Soc., 415pp, 2001
\bibitem[FS1]{fs1} T. Foertsch \& V. Schroeder, {\it Hyperbolic Rank of Products}, to appear in the Proc. Amer. Math. Soc.
\bibitem[FS2]{fs2} T.Foertsch \& V.Schroeder, {\it Products of hyperbolic metric spaces}, to appear in the Geom. Dedicata
\bibitem[G]{g} M.Gromov, {\it Hyperbolic groups}, Essays in group theory (S.M.Gersten, ed), Springer Verlag, MSRI Publ. 8 (1987),
75-263
\bibitem[M]{m} C.R.F. Maunder {\it Algebraic Topology}, Cambridge University Press, 1980
\end{thebibliography}
\end{document}